\documentclass[10pt]{amsart}
\usepackage{amssymb,amsmath,amscd,amsfonts,amsthm,mathrsfs}
\usepackage{verbatim,paralist}
\usepackage{enumerate}
\RequirePackage{ifthen}
\usepackage{tikz}
\usetikzlibrary{arrows}
\usepackage[english]{babel}
\usepackage{latexsym}
\usepackage{amsmath}
\usepackage{amsfonts}
\usepackage{amssymb}
\usepackage{latexsym}
\usepackage{color}

\newcommand{\Ac}{\mathcal{A}}
\newcommand{\Cc}{\mathcal{C}}
\newcommand{\Dc}{\mathcal{D}}
\newcommand{\Ec}{\mathcal{E}}
\newcommand{\Fc}{\mathcal{F}}
\newcommand{\Gc}{\mathcal{G}}
\newcommand{\Nc}{\mathcal{N}}
\newcommand{\Oc}{\mathcal{O}}
\newcommand{\Sc}{\mathcal{S}}
\newcommand{\Vc}{\mathcal{V}}
\newcommand{\off}{{\rm off}}

\newcommand{\N}{\mathbb{N}}
\newcommand{\C}{\mathbb{C}}
\newcommand{\Z}{\mathbb{Z}}
\newcommand{\id}{\rm id}
\newcommand{\supp}{{\rm supp}\, }

\newtheorem{theorem}{Theorem}[section]%{\bf}{\it }
\newtheorem{proposition}[theorem]{Proposition}%{\bf}{\it }
\newtheorem{remark}[theorem]{Remark}%{\bf}{\it }
\newtheorem{example}[theorem]{Example}%{\bf}{\it }
\newtheorem{definition}[theorem]{Definition}%{\bf}{\it }
\newtheorem{corollary}[theorem]{Corollary}%{\bf}{\it }
\title{The adjacency matrix and the discrete Laplacian acting on forms }
\author{Hatem Baloudi}
\address{Hatem Baloudi and Aref Jeribi, Facult\'e des sciences de
  sfax, Route de soukra km 3,5, B.P. 1171, 3000 Sfax, Tunisia} 
\author{Sylvain Gol\'enia}
\address{Sylvain Gol\'enia, Institut de Math\'ematiques, 351 cours de la
  lib\'eration, F 33405 Talence cedex, France}
\author{Aref Jeribi}


\begin{document}

\maketitle
\begin{abstract}
We study the relationship between the adjacency matrix and the
discrete Laplacian acting on $1$-forms.  We also prove that if the
adjacency matrix is bounded from below it is not necessarily
essentially self-adjoint. We discuss the question of essential
self-adjointness and the notion of completeness.
\end{abstract}

\section{Introduction}

Given a simple graph, the adjacency matrix $\Ac$ is unbounded from above if
and only if the degree of the graph is unbounded. It is more
complicated to characterize the fact that it is unbounded from
below. In \cite{1} one gives an optimal and necessary. The examples
provided in \cite{1} satisfy~:
\begin{enumerate}[1)]
\item $\Ac$ is not bounded from above,
\item $\Ac$ is bounded from below,
\item $\Ac$ is essentially self-adjoint.
\end{enumerate}
In \cite{1}, one asks the question. Can we have $1$ and $2$ without
having $3$? In Corollary \ref{c:az} we answer positively to this
question. Moreover, relying on different ideas than in \cite{1}, we
provide a large class of graphs satisfying $1$ and $2$. Our approach is
based on establishing a link between the adjacency matrices and the
discrete Laplacian $\Delta_1$ acting on $1$-forms.

The article is organized as follows: In Section 2, we
present general definitions about graphs. Then we define five
different types of discrete operator associated to a locally finite
graphs (the adjacency matrix, two discrete Laplacians acting on
vertices, two discrete Laplacians acting on edges). The link between
these definitions is made in Section \ref{edge}.
In Section \ref{s:first}, we start with the question of unboundedness
and give a upper bound on the infimum of the essential spectrum.
In Section \ref{s:esssa1} we provide a counter example to essential
 self-adjointness for the Laplacian acting on $1$-forms.
In Section \ref{s:complete} we present the notion of $\chi$-completeness
which was introduced in \cite{25} and develop it through optimal
examples. Then in Section \ref{s:other} we confront this criterion to
other type of approaches. Finally in Section \ref{s:adj} we provide a
new criterion of essential self-adjointness for the adjacency matrix.

\textbf{Acknowledgment}.
SG was partially supported by the ANR project
GeRaSic (ANR-13-BS01-0007-01) and SQFT (ANR-12-JS01-0008-01). HB
enjoyed the hospitality of Bordeaux University
when this work started. We would like to thank Michel Bonnefont for
useful discussions.

\section{Graph structures}

\subsection{Generalities about graphs}

We recall some standard definitions of graph theory. A \emph{graph} is
a triple
$\Gc:=(\mathcal{E},\mathcal{V}, m)$, where $\mathcal{V}$ is countable
set (the \emph{vertices}), $\mathcal{E}:\mathcal{V} \times
\mathcal{V}\rightarrow
\mathbb{R}_{+}$ (the \emph{edges}) is symmetric, and $m:\Vc\to (0, \infty)$ is a weight. When $m=1$, we simply write $\Gc=(\Ec, \Vc)$. We say that $\Gc$ is \emph{simple} if $m=1$ and $\Ec:\Vc\times \Vc\to \{0,1\}$.

Given  $x,y\in \mathcal{V}$, we say that $x$ and $y$ are
\emph{neighbors} if $\mathcal{E}(x,y)>0$. We denote this
relationship by $x\sim y$ and the set of neighbors of $x$ by
$\Nc_\Gc(x)$.
We say there is a \emph{loop} in $x\in\mathcal{V}$
if $\mathcal{E}(x,x)>0$.   A graph is \emph{connected} if for all
$x,y\in\mathcal{V}$, there exists a path $\gamma$ joining $x$ and
$y$. Here, $\gamma$ is a sequence
$x_{0},~x_{1},...,x_{n}\in\mathcal{V}$ such that $x=x_{0},~y=x_{n}$,
and $x_{j}\sim x_{j+1}$ for all $0\leq j \leq n-1$. In this case, we
endow $\mathcal{V}$ with the metric $\rho_{\mathcal{G}}$ defined by
\begin{align}\label{e:dist}
\rho_{\mathcal{G}}(x,y):=\inf
\{|\gamma|,~\gamma~\mbox{is~a~path~joining~$x$~and~$y$}\}.
\end{align}
A $x-$triangle given by $(x,y,z,x)$ is different form the one given by
  $(x,z,y,x)$.

  A graph $\Gc$ is \emph{locally finite} if $d_{\Gc}(x)$ is finite for
all $x\in\mathcal{V}$. In the sequel, we assume that:
\begin{center}
{\bf All graphs are locally finite, connected with no loops.}
\end{center}

 A graph $\Gc=(\Ec, \Vc, m)$ is  \emph{bi-partite}
  if there are $\Vc_1$ and $\Vc_2$ such that $\Vc_1\cap
  \Vc_2=\emptyset$, $\Vc_1\cup \Vc_2= \Vc$, and such that $\Ec(x,y)=0$
  for all $(x,y)\in \Vc_1^{2} \cup  \Vc_2^{2}$.

\begin{definition}\label{d:line} Given a graph $\Gc=(\Ec, \Vc, m)$, we
  define the
  \emph{line graph} of $\Gc$ the graph
  $\widetilde{\Gc}:=(\widetilde{\mathcal{E}},\widetilde{\mathcal{V}})$
  where $\widetilde{\mathcal{V}}$ is the set of edges of $\Gc$ and
\begin{align*}
\widetilde{\mathcal{E}}\Big((x_{0},y_{0}),(x,y)\Big)=\frac{\sqrt{\mathcal{E}(x_{0},y_{0})\mathcal{E}(x,y)}}{m(x)}1_{x=x_0}+
\frac{\sqrt{\mathcal{E}(x_{0},y_{0})\mathcal{E}(x,y)}}{m(y)}1_{y=y_0},
\end{align*}
if $(x_0,y_0)\neq(x,y)$ and $0$ otherwise.
\end{definition}
The choice of the weight $\tilde \Ec$ seems a bit arbitrary but is
motivated by Proposition \ref{p:uni}. We stress that $\tilde
m=1$. Others choices can be made, e.g., \cite{EvLa}. We refer to
\cite{Ha}[Chapter 8] for a presentation of this class in the setting
of simple graphs.

\begin{align*}
\begin{tikzpicture}[ scale=0.4]
\fill[color=black](-6, -6.49)circle(.5mm);
\fill[color=black](-5, -4.49)circle(.5mm);
\fill[color=black](-2.25, -3.70)circle(.5mm);
\fill[color=black](-0.10, -4.80)circle(.5mm);
\fill[color=black](-0.01, -8.49)circle(.5mm);
\fill[color=black](-2.41, -9.49)circle(.5mm);
\fill[color=black](-5.21, -8.90)circle(.5mm);
\fill[color=black](6, -6.49)circle(.5mm);
\fill[color=black](5, -4.49)circle(.5mm);
\fill[color=black](6.25, -7.70)circle(.5mm);
\fill[color=black](6.10, -8.80)circle(.5mm);
\fill[color=black](-3, -6.49)circle(.5mm);
\fill[color=black](2, -6.49)circle(.5mm);
\draw[-](-3, -6.49)--(2, -6.49);
\draw[-](-3, -6.49)--(-6, -6.49);
\draw[-](-3, -6.49)--(-5, -4.49);
\draw[-](-3, -6.49)--(-2.25, -3.70);
\draw[-](-3, -6.49)--(-0.10, -4.80);
\draw[-](-3, -6.49)--(-0.01, -8.49);
\draw[-](-3, -6.49)--(-2.41, -9.49);
\draw[-](-3, -6.49)--(-5.21, -8.90);
\draw[-](2, -6.49)--(6.10, -8.80);
\draw[-](2, -6.49)--(6.25, -7.70);
\draw[-](2, -6.49)--(5, -4.49);
\draw[-](2, -6.49)--(6, -6.49);
\path(-0.50, -10.95) node {A  bipartite graph $\Gc$};
\end{tikzpicture}
\quad \quad
\begin{tikzpicture}[ scale=0.4]
\fill[color=black](-6, -6.49)circle(.5mm);
\fill[color=black](-5, -4.49)circle(.5mm);
\fill[color=black](-2.25, -3.70)circle(.5mm);
\fill[color=black](-0.10, -4.80)circle(.5mm);
\fill[color=black](0.60, -6.49)circle(.5mm);
\fill[color=black](-0.01, -8.49)circle(.5mm);
\fill[color=black](-2.41, -9.49)circle(.5mm);
\fill[color=black](-5.21, -8.90)circle(.5mm);
\draw[-](-6, -6.49)--(-5, -4.49);
\draw[-](-5, -4.49)--(-2.25, -3.70);
\draw[-](-2.25, -3.70)--(-0.10, -4.80);
\draw[-](-0.10, -4.80)--(0.60, -6.49);
\draw[-](0.60, -6.49)--(-0.01, -8.49);
\draw[-](-0.01, -8.49)--(-2.41, -9.49);
\draw[-](-2.41, -9.49)--(-5.21, -8.90);
\draw[-](-5.21, -8.90)--(-6, -6.49);
\draw[-](-5, -4.49)--(-0.10, -4.80);
\draw[-](-0.10, -4.80)--(-0.01, -8.49);
\draw[-](-0.01, -8.49)--(-5.21, -8.90);
\draw[-](-5.21, -8.90)--(-6, -6.49);
\draw[-](-6, -6.49)--(-5, -4.49);
\draw[-](-5, -4.49)--(-5.21, -8.90);
\draw[-](-6, -6.49)--(-2.25, -3.70);
\draw[-](-6, -6.49)--(-0.10, -4.80);
\draw[-](-6, -6.49)--(0.60, -6.49);
\draw[-](-6, -6.49)--(-0.01, -8.49);
\draw[-](-6, -6.49)--(-2.41, -9.49);
\draw[-](-5, -4.49)--(0.60, -6.49);
\draw[-](-5, -4.49)--(-0.01, -8.49);
\draw[-](-5, -4.49)--(-2.41, -9.49);
\draw[-](-2.25, -3.70)--(0.60, -6.49);
\draw[-](-2.25, -3.70)--(-0.01, -8.49);
\draw[-](-2.25, -3.70)--(-2.41, -9.49);
\draw[-](-2.25, -3.70)--(-5.21, -8.90);
\draw[-](-0.10, -4.80)--(-2.41, -9.49);
\draw[-](-0.10, -4.80)--(-5.21, -8.90);
\draw[-](0.60, -6.49)--(-2.41, -9.49);
\draw[-](0.60, -6.49)--(-5.21, -8.90);
\fill[color=black](6, -6.49)circle(.5mm);
\fill[color=black](5, -4.49)circle(.5mm);
\fill[color=black](6.25, -7.70)circle(.5mm);
\fill[color=black](6.10, -8.80)circle(.5mm);
\draw[-](0.60, -6.49)--(6, -6.49);
\draw[-](0.60, -6.49)--(5, -4.49);
\draw[-](0.60, -6.49)--(6.25, -7.70);
\draw[-](0.60, -6.49)--(6.10, -8.80);
\draw[-](6, -6.49)--(5, -4.49);
\draw[-](6, -6.49)--(6.25, -7.70);
\draw[-](6.25, -7.70)--(6.10, -8.80);
\draw[-](5, -4.49)--(6.10, -8.80);
\draw[-](5, -4.49)--(6.25, -7.70);
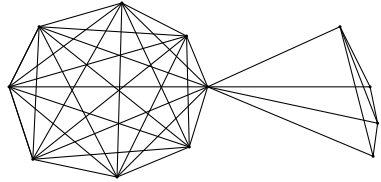
\path(-0.50, -10.95) node {The line-graph of $\Gc$};
\end{tikzpicture}
\end{align*}

\begin{align*}
\begin{tikzpicture}[scale=0.6]
\fill[color=black](-10.5, -6.49)circle(.5mm);
\fill[color=black](-14, -6.49)circle(.5mm);
\fill[color=black](-13, -4.49)circle(.5mm);
\fill[color=black](-10.25, -3.70)circle(.5mm);
\fill[color=black](-7.80, -4.80)circle(.5mm);
\fill[color=black](-7.10, -6.49)circle(.5mm);
\fill[color=black](-8.01, -8.49)circle(.5mm);
\fill[color=black](-10.41, -9.49)circle(.5mm);
\fill[color=black](-13.21, -8.90)circle(.5mm);
\draw[-](-10.5, -6.49)--(-14, -6.49);
\draw[-](-10.5, -6.49)--(-13, -4.49);
\draw[-](-10.5, -6.49)--(-10.25, -3.70);
\draw[-](-10.5, -6.49)--(-7.80, -4.80);
\draw[-](-10.5, -6.49)--(-7.10, -6.49);
\draw[-](-10.5, -6.49)--(-8.01, -8.49);
\draw[-](-10.5, -6.49)--(-10.41, -9.49);
\draw[-](-10.5, -6.49)--(-13.21, -8.90);
\path (-10.41, -10.95) node {A graph $\mathcal{K}_{1,8}$};
\fill[color=black](-6, -6.49)circle(.5mm);
\fill[color=black](-5, -4.49)circle(.5mm);
\fill[color=black](-2.25, -3.70)circle(.5mm);
\fill[color=black](-0.10, -4.80)circle(.5mm);
\fill[color=black](0.60, -6.49)circle(.5mm);
\fill[color=black](-0.01, -8.49)circle(.5mm);
\fill[color=black](-2.41, -9.49)circle(.5mm);
\fill[color=black](-5.21, -8.90)circle(.5mm);
\path(-2.50, -10.95) node { $\widetilde{\mathcal{K}_{1,8}}$ is the complete graph $K_{8}$};
\draw[-](-6, -6.49)--(-5, -4.49);
\draw[-](-5, -4.49)--(-2.25, -3.70);
\draw[-](-2.25, -3.70)--(-0.10, -4.80);
\draw[-](-0.10, -4.80)--(0.60, -6.49);
\draw[-](0.60, -6.49)--(-0.01, -8.49);
\draw[-](-0.01, -8.49)--(-2.41, -9.49);
\draw[-](-2.41, -9.49)--(-5.21, -8.90);
\draw[-](-5.21, -8.90)--(-6, -6.49);
\draw[-](-5, -4.49)--(-0.10, -4.80);
\draw[-](-0.10, -4.80)--(-0.01, -8.49);
\draw[-](-0.01, -8.49)--(-5.21, -8.90);
\draw[-](-5.21, -8.90)--(-6, -6.49);
\draw[-](-6, -6.49)--(-5, -4.49);
\draw[-](-5, -4.49)--(-5.21, -8.90);
\draw[-](-6, -6.49)--(-2.25, -3.70);
\draw[-](-6, -6.49)--(-0.10, -4.80);
\draw[-](-6, -6.49)--(0.60, -6.49);
\draw[-](-6, -6.49)--(-0.01, -8.49);
\draw[-](-6, -6.49)--(-2.41, -9.49);
\draw[-](-5, -4.49)--(0.60, -6.49);
\draw[-](-5, -4.49)--(-0.01, -8.49);
\draw[-](-5, -4.49)--(-2.41, -9.49);
\draw[-](-2.25, -3.70)--(0.60, -6.49);
\draw[-](-2.25, -3.70)--(-0.01, -8.49);
\draw[-](-2.25, -3.70)--(-2.41, -9.49);
\draw[-](-2.25, -3.70)--(-5.21, -8.90);
\draw[-](-0.10, -4.80)--(-2.41, -9.49);
\draw[-](-0.10, -4.80)--(-5.21, -8.90);
\draw[-](0.60, -6.49)--(-2.41, -9.49);
\draw[-](0.60, -6.49)--(-5.21, -8.90);
\end{tikzpicture}
\end{align*}

As we see on this example, the line-graph is obtained by gluing
together some complete graphs. This is not true in general.
The classification of the line-graphs of a finite graph is well-known, see
\cite{Ha}[Theorem 8.4].

We present the relationship between bipartite graphs and line graphs.
We recall that $C_{n}$ denotes the $n$-cycle graph, i.e., $\Vc:=\Z / n\Z$,
where $\Ec(x,y)=1$ if and only if $|x-y|=1$.

\begin{proposition} Let $\Gc=(\mathcal{E},\mathcal{V})$ be a bipartite
  graph. Then the two following assertions are equivalent:
\begin{enumerate}[1)]
\item $\Gc\simeq \widetilde{\Gc}$,
\item $\Gc\in \{\mathbb{Z},~\mathbb{N},~C_{2n}:~n\in\mathbb{N}\}$.
\end{enumerate}
\end{proposition}
\proof
 It clear that $(2)\Rightarrow (1)$. Now, suppose that $\Gc\simeq
 \widetilde{\Gc}$. Then $\widetilde{\Gc}$ not contain $x-$triangles. If
 $d_{\Gc}(x)\geq 3$, for all $x\in\mathcal{V}$, then there exist a
 $x-$triangle of $\widetilde{\Gc}$. So, $d_{\Gc}(x)\leq 2~\forall
 x\in\mathcal{V}$. This proof is complete. \qed

\subsection{The adjacency matrix}\label{s:adjacency}
Set $\Gc=(\Ec, \Vc,m)$ a weighted graph. We define the set of $0-$cochains
$\mathcal{C}(\Vc):= \{f:\Vc\to \C\}$. We denote by $\Cc^c(\Vc)$ the $0$-cochains with finite support in $\Vc$. We associate a Hilbert space to $\Vc$.
\begin{align*}
\ell^{2}(\Vc):=\left\{f\in\mathcal{C}(\Vc)\mbox{~such~that~}\|
f \|^{2}:\displaystyle=\sum_{x\in
  \mathcal{V}} m(x)| f(x)|^{2}<\infty\right\}.
\end{align*}
The associated scalar product is given by
\begin{align*}
 \langle f,g\rangle:=\sum_{x\in\mathcal{V}} m(x)\overline{f(x)}g(x),\mbox{~for~}f,g\in\ell^{2}(\mathcal{V}).
\end{align*}
We define the \emph{adjacency matrix}:
\begin{align*}
\mathcal{A}_{\Gc}(f)(x):= \frac{1}{m(x)}\sum_{y\in\Vc}\mathcal{E}(x,y)f(y)
\end{align*}
for $f\in   \mathcal{C}^{c}(\Vc)$. It is symmetric and thus
closable. We denote its closure by the same symbol.
We denote the \emph{degree} by
\[\deg_\Gc(x):=\frac{1}{m(x)}\sum_{y\in \Vc}\Ec(x,y).\]
When $m=1$, we
have that $\Ac_\Gc$ is unbounded  if and only it is unbounded from
above and if and only if the degree is unbounded. The
fact that $\Ac_\Gc$ is unbounded from below or not is a delicate
question, see \cite{1}. We refer to \cite{9} for the description of some
self-adjoint extensions.

\section{The laplacians acting on edges}

\subsection{The symmetric and anti-symmetric spaces}
In the previous section we defined a Hilbert structure on the set of
vertices. We can endow two different Hilbert structures on the set of
edges. Both have their own interest.

Given a weighted graph $\Gc= (\Ec, \Vc,m)$, the set of $1-$cochains (or
$1$-forms) is given by:
\[ \mathcal{C}_{{\rm anti}}(\Ec):=  \{f:\Vc\times \Vc \to \C,\,
f(x,y)=-f(y,x) \mbox{ for all } x,y\in \Vc\},\]
where anti stands for anti-symmetric. This corresponds to fermionic
statistics. Concerning bosonic statistics, we define:
\[ \mathcal{C}_{{\rm sym}}(\Ec):=  \{f:\Vc\times \Vc \to \C,\,
f(x,y)=f(y,x) \mbox{ for all } x,y\in \Vc\}.\]
The sets of functions with finite support are denoted by
$\mathcal{C}_{{\rm anti}}^{c}(\Ec)$ and
$\mathcal{C}_{{\rm sym}}^{c}(\Ec)$,  respectively.

We turn to the Hilbert structures.
\begin{align*}
\ell^{2}_{{\rm anti}}(\mathcal{E}):=\left\{f\in\mathcal{C}_{{\rm anti}}(\Ec)\mbox{~such~that~}\|
f \|^{2}:\displaystyle=\frac{1}{2}\sum_{x,y\in
  \mathcal{V}}\mathcal{E}(x,y) | f(x,y)|^{2}<\infty\right\}
\end{align*}
and
\begin{align*}
\ell^{2}_{{\rm sym}}(\mathcal{E}):=\left\{f\in\mathcal{C}_{{\rm sym}}(\Ec)\mbox{~such~that~}\|
f \|^{2}:\displaystyle=\frac{1}{2}\sum_{x,y\in
  \mathcal{V}}\mathcal{E}(x,y) | f(x,y)|^{2}<\infty\right\}.
\end{align*}
The associated scalar product is given by
\begin{align*}
 \langle
 f,g\rangle:=\frac{1}{2}\sum_{x,y\in\mathcal{V}}\mathcal{E}(x,y)\overline{f(x,y)}g(x,y),
\end{align*}
when $f$ and $g$ are both in $\ell^{2}_{{\rm anti}}(\mathcal{E})$ or
in $\ell^{2}_{{\rm sym}}(\mathcal{E})$.

We start by defining operator in the anti-symmetric case.   The
\emph{difference operator} is defined as
\begin{align*}
d:=d_{\rm anti}:\mathcal{C}^{c}(\Vc)\longrightarrow \mathcal{C}_{{\rm
    anti}}^{c}(\Ec),~d(f)(x,y):=f(y)-f(x). \end{align*}
The \emph{coboundary operator} is the formal adjoint of $d$. We set:
\begin{align*}
d^*:=
%d^{*}_{\rm sym}
d^{*}_{\rm anti}:\mathcal{C}_{{\rm anti}}^{c}(\Ec)\longrightarrow
\mathcal{C}^{c}(\Vc),~d^{*}(f)(x):=\frac{1}{m(x)}\sum_{y\in\mathcal{V}}\mathcal{E}(y,x)f(y,x).
\end{align*}
We denote by the same symbols their closures.  The Gau{\ss}-Bonnet operator is defined on
  $\mathcal{C}^{c}(\Vc)\oplus \mathcal{C}_{{\rm anti}}^{c}(\Ec)$ by
\begin{align*}
D:=d+d^{*}\cong \left(
\begin{array}{cc}
      0 & d^{*}\\
      d & 0\\
\end{array}
\right).\end{align*}
This operator is motivated by geometry, see \cite{25}. It is also a Dirac-like
operator.
 The associated Laplacian is defined as $\Delta:=D^{2}=
 \Delta_{0}\oplus\Delta_{1}$ where $\Delta_{0}$ is the  discrete
 Laplacian acting on $0-$forms given by
\begin{equation}
\Delta_{0}(f)(x):=(d^*d)
(f)(x)=\frac{1}{m(x)}\sum_{y}\mathcal{E}(x,y)(f(x)-f(y)),
\end{equation}
with $f\in \mathcal{C}^{c}(\Vc)$,
and where the discrete Laplacian acting on $1-$forms is given by
\begin{align}\nonumber
\Delta_{1}(f)(x,y)&:=(dd^*)(f)(x,y)
\\ \label{144}
&=
\frac{1}{m(x)}\sum_{z\in\mathcal{V}}\mathcal{E}(x,z)f(x,z)+\frac{1}{m(y)}\sum_{z\in\mathcal{V}}\mathcal{E}(z,y)f(z,y),
\end{align}
with $f\in  \mathcal{C}_{\rm anti}^{c}(\Ec)$. Both operator are
symmetric and thus closable.  We denote the closure by $\Delta_0$
(resp.\ $\Delta_1$), its domain by
$\mathcal{D}(\Delta_0)$ (resp.\ $\mathcal{D}(\Delta_1)$), and its
adjoint by $\Delta_{0}^{*}$ (resp.\ $\Delta_{1}^{*}$).

 In \cite{16}, one proves that $\Delta_0$ is essentially self-adjoint
 on $\Cc^c(\Vc)$ when the graph is simple. The literature on the
 subject is large and refer to \cite{Gol2} for historical purposes.
 In \cite{25}, they prove
 that $D$ is essentially
 self-adjoint on $\Cc^c(\Vc)\oplus \Cc_{\rm anti}^c(\Ec)$ if and only if
$\Delta_0\oplus\Delta_1$ is essentially self-adjoint on $\Cc^c(\Vc)\oplus
  \Cc_{\rm anti}^c(\Ec)$. Moreover
 they provide a criterion based on the completeness of a metric. In
 Section \ref{s:esssa1} we provide a counter example to essential
 self-adjointness for the Laplacian acting of $1$-forms. Then in
 Section \ref{s:complete} we discuss the notion of $\chi$-completeness
 that was introduced in \cite{25}. This notion is related to the one of
 intrinsic metric, e.g., \cite{X.} and references therein.

We turn to the symmetric choice.   We set~:
\begin{align*}
d_{\rm sym}:\mathcal{C}^{c}(\Vc)\longrightarrow \mathcal{C}_{{\rm
    sym}}^{c}(\Ec),~d_{\rm sym}(f)(x,y):=f(y)+f(x).
\end{align*}
It is interesting to note that $d_{\rm sym}$ is in fact the \emph{incidence
matrix} of $\Gc$. Indeed
\[d_{\rm sym}(\delta_{x_{0}})(x,y)= \left\{\begin{array}{ll}
1, & \mbox{if } x_0\in \{x,y\},
\\
0, & \mbox{otherwise}.
\end{array}\right. \]
The formal adjoint of the incidence matrix is given by:
\begin{align*}
d^{*}_{\rm sym}:\mathcal{C}_{{\rm sym}}^{c}(\Ec)\longrightarrow
\mathcal{C}^{c}(\Vc),~d^{*}_{\rm
  sym}(f)(x):=\frac{1}{m(x)}\sum_{y\in\mathcal{V}}\mathcal{E}(y,x)f(y,x).
\end{align*}
A direct computation gives that~:
\[\Delta_{0, \rm sym} f(x):= (d^{*}_{\rm sym}d_{\rm sym}) f(x)= \frac{1}{m(x)}\sum_{y} \Ec(x,y)(f(x)+f(y)),\]
for all $f\in\Cc_{}^c(\Vc)$ and
\begin{align}\nonumber
\Delta_{1, \rm sym} f(x,y)&:= (d_{\rm sym}d^{*}_{\rm sym}) f(x,y)
\\ \label{e:adjcentral}
&= \frac{1}{m(x)}\sum_{z\in \Vc} \Ec(x,z) f(x,z) +
\frac{1}{m(y)}\sum_{z\in \Vc} \Ec(z,y) f(z,y),
\end{align}
for all $f\in\Cc_{\rm sym}^c(\Ec)$.

We denote by the same symbol the closure of these operators. We turn to this central observation.

\begin{proposition}\label{p:uni}
Set  $\Gc=(\Ec, \Vc, m)$. The operator $\Delta_{1,\rm sym}$ is
unitarily equivalent to
\[\Ac_{\tilde \Gc} + V(\cdot, \cdot),\]
where
\[V(x,y):=\Ec(x, y)\left(\frac{1}{m(x)}+ \frac{1}{m(y)}\right)\]
and $\tilde \Gc$ denotes the line-graph of $\Gc$, see Definition \ref{d:line}.
\end{proposition}

It is important to note that $\ell^2_{\rm sym}(\Ec)$ is a weighted
space, whereas $\ell^2(\tilde \Vc)$ is not one. This proposition is
part of the folklore in the case of finite and unweighted graphs.

\proof We set $U: \ell^2_{\rm sym}(\Ec) \to \ell^2(\tilde \Vc, 1)$ as
being the unitary transformation given by $Uf(x,y):= \sqrt{\Ec(x,y)}
f(x,y)$. We have $(U \Delta_{1, \rm sym} U^{-1}) f(x,y)=$
\begin{align*}
\sum_{z\in \Vc} \left(\frac{\sqrt{\Ec(x,y)\Ec(x,z)}}{m(x)} f(x,z) + \frac{\sqrt{\Ec(x,y)\Ec(z,y)}}{m(y)} f(z,y)\right),
\end{align*}
for all $f\in \Cc^c(\widetilde{\Vc})$. Remembering the definition
of the weights for a line-graph, we obtain the result. \qed

\begin{corollary}
Let $\Gc=(\Ec, \Vc, m)$ be a weighted graph. If
\[c:=\sup_{x,y\in \Vc} \Ec(x, y)\left(\frac{1}{m(x)}+ \frac{1}{m(y)}\right) <\infty\]
then $\Ac_{\tilde \Gc}\geq -c$. In particular, $\Ac_{\tilde \Gc}\geq
-2$ when $\Gc$ is simple.
\end{corollary}

In the setting of simple graphs, in \cite{1}, one gives a necessary
condition for the semi-boundedness of the adjacency matrix. The
\emph{lower local complexity} of a graph $\Gc$ is defined by
\begin{align*}
C_{\rm loc}(\Gc):=\inf \bigcap \left\{\overline{\Big\{\displaystyle
  \frac{N_{\Gc}(x)}{d_{\Gc}^{2}(x)},~x\in\mathcal{V}~\textup{and}~d_{\Gc}(x)\geq
  n\Big\}},~n\in\mathbb{N}\right\},
\end{align*}
  where $N_{\Gc}(x):=\sharp \{x-\textup{triangles}\}$. The
  \emph{sub-lower local complexity} of a graph $\Gc$
  is defined by
\begin{align*}
C_{\rm loc}^{\rm sub}(\Gc)=\inf_{\{\Gc^{'}\subset \Gc,~\sup
  d_{\Gc^{'}}=\infty\}}C_{\rm loc}(\Gc^{'}),
\end{align*}
where $\Gc'$ is a sub-graph of $\Gc$, i.e., $\Gc'$ is a graph whose
vertex set is a subset of that of $\Gc$, and whose adjacency relation
is a subset of that of $\Gc$ restricted to this subset.

We recall \cite[Theorem 1.1]{1}:

\begin{theorem}
Let $\Gc=(\Ec,\Vc)$ be a locally finite graph  of unbounded class. Let
$\hat \Ac_{\Gc}$ be a self-adjoint realization of the $\Ac_{\Gc}$. Suppose that
$\Ec$ is bounded. Then, one has:
\begin{enumerate}[1)]
\item $\hat \Ac_{\Gc}$ is unbounded from above.
\item If  $C_{\rm loc}^{\rm sub}(G)=0$, then $\hat \Ac_{\Gc}$
  is unbounded from
  below.
\item For all $\varepsilon>0$, there is a connected simple graph $\Gc$
  such that $C_{\rm  loc}(G)\in (0, \varepsilon)$, $\Ac_{\Gc}$
  is essentially
  self-adjoint on $\Cc^c(\Vc)$  and is bounded from below.
\end{enumerate}
\end{theorem}

The line-graph do not lead to the optimality that is obtained in third
point. Indeed, by repeating the proof of \cite[Proposition 3.2]{1} we
obtain that given  $\Gc=(\Ec, \Vc)$ be a locally finite simple
bipartite graph, then
\[C_{\rm loc}^{\rm
  sub}(\tilde\Gc)>\frac{1}{10}.\]

\subsection{Breaking symmetry}\label{edge}
We have that \eqref{144} and \eqref{e:adjcentral} have the same
expression. However they do not act on the same spaces. When $\Gc$ is
bi-partite, we shall prove that the two operators are unitarily equivalent.

We now fix an orientation of a graph $\Gc=(\mathcal{E},\mathcal{V}, m)$. For each edge there are two possible choices. For $x,y\in\mathcal{V}$ such that
$\mathcal{E}(x,y)\neq 0$, we write $x\rightarrow y$ or $y\rightarrow
x$, following the choice of the orientation.

Let $U:\ell^{2}_{\rm anti}(\mathcal{E})\longrightarrow \ell_{\rm sym}^{2}(\mathcal{E})$ be the unitary map given by
\begin{align*}
(Uf)(x,y)={\rm sign}(x,y)f(x,y), \quad \mbox{where }
{\rm sign}(x,y):=\left\{
\begin{array}{rl}
1, &  \textup{if}~x\rightarrow y ,\\
-1, & \textup{if}~y\rightarrow x.
\end{array}
\right.
\end{align*}
Note that $(U^{-1}f)(x,y)={\rm sign}(x,y)f(x,y)$. Therefore, for $x_{0}, y_{0}\in\mathcal{V}$ and $f\in\ell^{2}_{\rm sym}(\mathcal{E})$ we have
\begin{align}\nonumber
U\Delta_{1}U^{-1}f(x_{0},y_{0})=
&\, {\rm sign}(x_{0},y_{0})
\Big(\frac{1}{m(y_{0})}\sum_{x\in \Vc}\mathcal{E}(x,y_{0}){\rm sign}(x,y_{0})f(x,y_{0})
\\ \label{e:choosesign}
& +\frac{1}{m(x_{0})}\sum_{y\in \Vc}\mathcal{E}(x_{0},y){\rm sign}(x_{0},y)f(x_{0},y)\Big).
\end{align}

In order to take advantage of this transformation, one has to choose
the orientation in a good way. We have~:

\begin{proposition}\label{p:unisym}
Assume that $\Gc=(\Ec, \Vc, m)$ is a bi-partite weighted space, then
$\Delta_1$ and $\Delta_{1,\rm sym}$ are unitarily equivalent. In particular, if we also have that $\Gc$ is simple, $\Delta_1$ is unitarily equivalent to $\Ac_{\tilde \Gc}+ 2 \id$.
\end{proposition}
\proof We consider the bi-partite decomposition $\{\Vc_1, \Vc_2\}$ of the graph $\Gc$. For $x\in \Vc_1$ and $y\in \Vc_2$, we set $x\to y$. From \eqref{e:choosesign}, we see that all the products of the signs are being $1$. This gives the first part. For the rest, apply Proposition \ref{p:uni}. \qed

Since the line-graph of $\Z$ is itself and it is also bi-partite, we can make a different choice of orientation.

\begin{example} We consider on $\Z$ the orientation $n\rightarrow n+1$, for all $n\in \Z$.
\begin{align*}
\begin{tikzpicture}[ scale=0.6]
\fill[color=black](-12.1, 0)circle(.5mm);
\fill[color=black](-14.1, 0)circle(.3mm);
\fill[color=black](-14.3, 0)circle(.3mm);
\fill[color=black](-14.5, 0)circle(.3mm);
\fill[color=black](-10.1, 0)circle(.5mm);
\fill[color=black](-8.1, 0)circle(.5mm);
\fill[color=black](-6.1, 0)circle(.5mm);
\fill[color=black](-4.1, 0)circle(.5mm);
\fill[color=black](-2.1, 0)circle(.3mm);
\fill[color=black](-2.3, 0)circle(.3mm);
\fill[color=black](-2.5, 0)circle(.3mm);
\draw[-triangle 60](-12.1, 0)--(-10.1, 0);
\draw[-triangle 60](-10.1, 0)--(-8.1, 0);
\draw[-triangle 60](-8.1, 0)--(-6.1, 0);
\draw[-triangle 60](-6.1, 0)--(-4.1, 0);
\end{tikzpicture}
\end{align*}
Note that $\tilde \Z \simeq \Z$. Moreover, we have:
\begin{align*} (U\Delta_{1}U^{-1})f(n,n+1)=2f(n,n+1)-\mathcal{A}_{\widetilde{\Gc}}f(n,n+1). \end{align*}
Therefore, $\Delta_{1}$ and $2I-\mathcal{A}_{\widetilde{\Gc}}$ are unitarily equivalent.
\end{example}

\section{On the spectral properties of the Laplacian  acting on
forms}\label{s:first}
We start with the question of unboundedness.

\begin{proposition}Let $\Gc=(\Ec,\Vc, m)$ be a weighted graph.
Then $\Delta_{1}$ is bounded if and only if $\Delta_{1, {\rm sym}}$ is
bounded if and only if $\Delta_{0}$ is bounded if and
only if
\[\displaystyle\sup_{x\in \Vc} \deg_\Gc(x)<\infty.\]
\end{proposition}

\proof The last equivalence is well-known, e.g., \cite{KL2}. On one side we have
\[0 \leq \langle f, \Delta_0 f\rangle \leq \langle f, 2\deg_\Gc(\cdot) f \rangle,\]
for all $f\in \Cc^c(\Vc)$ and on the other side $\langle \delta_x, \Delta_0 \delta_x\rangle= \deg_\Gc(x)$. The same statements hold true for $\Delta_{0, \rm sym}$.

We turn to $\Delta_1$. We denote by $\Fc$ the Friedrichs extension.
\begin{align*}
\Delta_1 \mbox{ is bounded } &\Leftrightarrow  \Delta_1^\Fc \mbox{ is bounded}
 \Leftrightarrow d^* \mbox{ is
  bounded }
  \\
  &\Leftrightarrow  d \mbox{ is
  bounded }
\Leftrightarrow  \Delta_0 \mbox{ is
  bounded}.
\end{align*}
The second equivalence comes from the fact that $\|(\Delta_1^\Fc)^{1/2}f\| = \|d^* f\|$, for all $f\in \Cc^c_{\rm anti}(\Ec)$ and by construction of the Friedrichs extension (e.g., \cite{RS}). The last one is of the same type.
The equivalence with $\Delta_{1, {\rm sym}}$ is
similar. \qed

Using the Persson's Lemma we give the following result:
\begin{proposition}
Let $\Gc=(\mathcal{E},\mathcal{V}, m)$ be an infinite and $\Delta_{1}^{\mathcal{F}}$ be a Friedrichs extension of $\Delta_{1}$. Then
\begin{align*}
\inf \sigma(\Delta_{1}^{\mathcal{F}})\leq\inf_{x,y\in \Vc, x\sim y}\left(
  \frac{1}{m(x)}+\frac{1}{m(y)}\right)\mathcal{E}(x,y)\end{align*}
\noindent and
\begin{align*}
\inf \sigma_{\rm ess}(\Delta^{\mathcal{F}}_{1})\leq \inf_{K\subset
  \mathcal{E}, K \rm finite}\inf_{(x,y)\in
  K^{c}, x\sim y} \left(\frac{1}{m(x)}+\frac{1}{m(y)}\right)\mathcal{E}(x,y).
\end{align*}
  In particular,  $\Delta_{1}^{\mathcal{F}}$
  is \emph{not} with compact   resolvent when $\Gc$ is simple and
  infinite.
\end{proposition}
 The last statement really differs from the situation case of
 $\Delta_0$. Indeed, in the case of a simple sparse graph (a planar graph for
 instance), one has that
 $\Dc(\Delta_0^{1/2})=\Dc\left(\deg_\Gc^{1/2}(\cdot)\right)$, see
 \cite{BGK} and also \cite{Gol2}. In particular we have that
 $\Delta_0$ is with compact resolvent if and only if
 $\lim_{|x|\to\infty} \deg_\Gc(x)=\infty$.

\proof
Let $x_{0},y_{0}\in\mathcal{V}$ such that $\mathcal{E}(x_{0},y_{0})\neq 0$. Then,
\begin{align*}
\left\langle\displaystyle\frac{\delta_{x_{0},y_{0}}-
\delta_{y_{0},x_{0}}}{\sqrt{\mathcal{E}(x_{0},y_{0})}},\Delta_{1}^{\mathcal{F}}\Big(\displaystyle\frac{\delta_{x_{0},y_{0}}-
\delta_{y_{0},x_{0}}}{\sqrt{\mathcal{E}(x_{0},y_{0})}}\Big)\right\rangle=
\left(\frac{1}{m(x_{0})}+\frac{1}{m(y_{0})}\right)\mathcal{E}(x_{0},y_{0}).\end{align*}
Recalling
\begin{align}\label{e:1}
\Big\|\frac{\delta_{x_{0},y_{0}}-\delta_{y_{0},x_{0}}}{\sqrt{\mathcal{E}(x_{0},y_{0})}}\Big\|=1.
\end{align}
 and the Persson's Lemma, e.g., \cite{KL2}[Proposition 18], the result follows. \qed

\begin{remark} The same result holds true for the operator $\Delta_{1, {\rm sym}}$ by
  changing the  $-$ into $+$ in \eqref{e:1}.
\end{remark}

\section{Essential self-adjointness}\label{s:esssa}
From now on, we concentrate on the analysis of $\Delta_{1}$. The
results for $\Delta_{1, {\rm sym}}$ are the same and the proofs need
only  minor changes.

\subsection{A counter example}\label{s:esssa1}
First it is important to notice that unlike $\Delta_0$, $\Delta_{1}$,
define as the closure of \eqref{144}, is not necessarily self-adjoint
on simple graphs.

Let $\Gc=(\mathcal{E},\mathcal{V})$ be radial simple tree. We denote the origin by $\nu$ and the spheres by $\Sc_n:=\{x\in \Vc,\, \rho_\Gc(\nu, x)=n\}$. Let $\off(n):=|\Sc_{n+1}|/|\Sc_n|$ be the \emph{offspring} of the $n$-th generation.

\begin{theorem}\label{az} Let $\Gc=(\mathcal{E},\mathcal{V})$ be a radial simple tree. Suppose that
\begin{align}\label{h}
n\mapsto \displaystyle \frac{\off^{2}(n)}{\off(n+1)}\in\ell^{1}(\mathbb{N}).
\end{align}
Then, $\Delta_{1}$  does not essentially self-adjoint on $\mathcal{C}_{\rm anti}^{c}(\Ec)$ and the deficiency indices are infinite.
\end{theorem}

\proof Set $f\in\ell^{2}_{\rm anti}(\mathcal{E})\setminus\{0\}$ such that
$f\in \ker(\Delta_{1}^{*}+{\rm i})$ and such that $f$ is constant on $\Sc_{n}\times \Sc_{n+1}$.
We denote the constant value by $C_{n}$.
So, we have the following equation
\begin{equation*}
(\off(n)+1+{\rm i})C_{n}-C_{n+1}\off(n+1)-C_{n-1}=0.
\end{equation*}
  Therefore,
\begin{align*}
\left\| f|_{\Sc_{n+1}\times \Sc_{n+2}}\right\|^2
&=| C_{n+1} |^{2}\prod_{i=0}^{n+1}\off(i)
\\
&\hspace{-0.7 cm}\leq 2\frac{| \off(n)+1-\rm i|^{2}}{\off^{2}(n+1)}\prod_{i=0}^{n+1}\off(i)| C_{n}|^{2}+2\frac{1}{\off^{2}(n+1)}\prod_{i=0}^{n+1}\off(i)| C_{n-1}|^{2}.
\\
&\hspace{-0.7 cm}= 2\frac{| \off(n)+1-\rm i|^{2}}{\off(n+1)}\|f|_{\Sc_{n}\times \Sc_{n+1}}\|^{2}+2\frac{\off(n)}{\off(n+1)}\|f|_{\Sc_{n-1}\times \Sc_{n}}\|^{2}.
\end{align*}
Since $\displaystyle \frac{\off^{2}(n)}{\off(n+1)}$ tends to $0$ as $n$ goes to infinity, we get by induction:
\[C:=\sup_{n\in\N} \|f|_{\Sc_{n}\times \Sc_{n+1}}\|^{2}<\infty.\]
Then, we have
\begin{align*}
\| f |_{\Sc_{n+1}\times \Sc_{n+2}}\|^2\leq 2C \left( \frac{| \off(n)+1-\rm i|^{2}}{\off(n+1)}+\frac{\off(n)}{\off(n+1)}
\right)
.\end{align*}
By $($\ref{h}$)$, we conclude that $f\in\ell^{2}(\mathcal{E})$ and
that $\dim \ker(\Delta_{1}^{*}+{\rm i})\geq 1$.    Using \cite[Theorem
X.36]{RS} we derive that $\Delta_{1}$ is not essentially self-adjoint.

Finally, since the deficiency indices are stable under bounded
perturbation, by cutting off the graph large enough, we reproduce the
same proof with an arbitrary large number of identical and disjoint
trees, we see that $\dim \ker(\Delta_{1}^{*}+{\rm i})=\infty$. \qed

\begin{remark} For the adjacency matrix, we are able to construct an
  example of simple graph where the deficiency indices are equal to
  $1$, see \cite{GoSc2}. The case of $\Delta_1$ is more complicated
  and we were not able to find an example where the deficiency indices
  are finite. We leave this question open.
\end{remark}

We can answer a problem which was left open in \cite{9}.

\begin{corollary}\label{c:az}
There exists a locally finite simple graph $\Gc=(\Ec, \Vc)$ such that $\Ac_\Gc$ is
bounded from below by $-2$ and such that $\Ac_\Gc$ is \emph{not} essentially
self-adjoint on $\Cc^c(\Vc)$.
\end{corollary}
\proof Combine Theorem \ref{az} with Proposition \ref{p:unisym}. \qed

\subsection{Completeness}\label{s:complete}

We first recall a criterion obtained in \cite{25}. They
introduce the following definition (see references therein for
historical purposes).

\begin{definition}
The graph  $\Gc=(\Ec, \Vc, m)$ is \emph{$\chi-$complete} if there
  exists a increasing sequence of finite set $(\Oc_{n})_{n}$ such that
  $\Vc=\cup_{n}\Oc_{n}$ and there exist related functions $\chi_{n}$
  satisfying the following three conditions:
\begin{enumerate}[1)]
\item $\chi_{n}\in\Cc^{c}(\Vc),~0\leq \chi_{n} \leq 1$,
\item $\chi_{n}(x)=1$ if $x\in \Oc_{n},$
\item $\exists C>0,~\forall n\in\mathbb{N},~x\in \Vc$, such that
\begin{align*}\displaystyle \frac{1}{m(x)}\sum_{y}\Ec(x,y)
  |\chi_{n}(x)-\chi_{n}(y)|^{2}\leq C.\end{align*}
\end{enumerate}
\end{definition}
%The third condition shares some similarities with \eqref{e:sim}.
The
result of \cite{25}, can be reformulate as follows:

\begin{theorem}\label{t:25}
Take $\Gc=(\Ec, \Vc, m)$ to be \emph{$\chi-$complete}. Then
$\Delta_1$ is essentially self-adjoint on $\Cc_{\rm anti}^c(\Ec)$.
\end{theorem}

\begin{remark} The converse is an open problem. The aim of the next sections
  is to tackle it.
\end{remark}

\subsubsection{A sharp example}

Theorem \ref{t:25} is abstract and may be difficult to check in
concrete settings. In \cite{25}, they provide only one example. Moreover, they
do not give a concrete way to prove that a graph is not
$\chi$-complete. To start-off we provide a sharp example.

\begin{proposition}\label{p:25sharp}
Let $\Gc=(\Ec, \Vc)$ be a locally finite simple tree,
endowed with an origin such that $\off(n):=\off(x)$ for
all $x\in \Sc_{n}$. Then the graph $\Gc$ is $\chi-$complete is and
only if
\[\sum_{n=1}^\infty \frac{1}{\sqrt{\off(n)}}=\infty.\]
\end{proposition}

\begin{example}\label{e:25sharp}
Set $\alpha>0$. Let $\Gc=(\Ec, \Vc)$ be a locally finite simple tree,
endowed with an origin such that $\off(x)=\lfloor n^\alpha\rfloor$ for
all $x\in \Sc_{n}$. Then the graph $\Gc$ is $\chi-$complete is and
only if $\alpha\leq 2$.
\end{example}

\proof
Suppose that the series converges and that $\Gc$ is $\chi-$complete. Let
$(\chi_n)_{n\in \N}$ be as in the definition. Using 3), we get:
\[|\chi_n(m)- \chi_n(m+1)|\leq \frac{\sqrt{C}}{\sqrt{\off(m)}}.\]
Moreover, by convergence of the series, there is $N\in \N$ such that
\[\sum_{k=N}^\infty\frac{1}{\sqrt{\lfloor n^\alpha
    \rfloor}}<\frac{1}{2\sqrt C}.\]
Then, by 2), there is $n_0\in \N$ such that $\chi_{n_0}(x)=1$ for all
$|x|\leq N$. Since $\chi_{n_0}$ is with finite support, there is $M\in
\N$ such that $\chi_{n_0}(x)=0$ for all $|x|\geq N+M$. Therefore,
\begin{align*}
1&=|\chi_{n_0}(N) -\chi_{n_0}(N+M)|
\\
&\leq |\chi_{n_0}(N) -\chi_{n_0}(N+1)| +\ldots + |\chi_{n_0}(N+M-1)
-\chi_{n_0}(N+M)|
\\
& \displaystyle \leq \sqrt{C}\sum_{n=N}^{N+M-1}\frac{1}{\sqrt{\lfloor
    n^\alpha\rfloor}}<\frac{1}{2}.
\end{align*}
Contradiction. Therefore if the series converges, then $\Gc$ is not
$\chi-$complete.

Suppose now that the series diverges. Set
\[\chi_n(x):=\left\{ \begin{array}{rl}
1, & \mbox{ if } |x|\leq n,
\\ \displaystyle
\max\left(0, 1- \sum_{k=n}^{|x|-1} \frac{1}{\sqrt{\off(k)}}\right),&
  \mbox{ if } |x|> n.
\end{array} \right.\]
Since the series diverges, $\chi_n$ is with
compact support and satisfies the definition of
$\chi$-completeness. \qed

\subsubsection{$1$-dimensional decomposition}

We now strengthen the previous example and follow ideas of \cite{BG}.

\begin{definition} A \emph{$1$-dimensional decomposition}
of the graph $\Gc:=(\Vc,\Ec)$ is
a family of finite sets $(S_n)_{n\geq 0}$ which forms a partition of
$\Vc$, that is $\Vc=\sqcup_{n\geq 0} S_n$, and such that for all $x\in S_n,
y\in S_m$,
\begin{eqnarray*}
 \Ec (x,y)>0 \implies |n-m|\leq 1.
\end{eqnarray*}
Given such a $1$-dimensional decomposition, we write $|x|:=n$ if $x\in
S_n$. We also write $B_N:=\cup_{0\leq i\leq N} S_i$.
\end{definition}
 We stress that
$S_0$ is not asked to be a single point.

\begin{example} Typical examples of such a $1$-dimensional decomposition
  are given by level sets of the graph  distance function to a finite
  set $S_0$ that is
\begin{equation}\label{e:1ddSn}
S_n:=\{x \in \Vc, \rho_\Gc(x,S_0)=n\},
\end{equation}
where the \emph{graph distance function} $\rho_\Gc$ is defined by
\eqref{e:dist}.
Note that for a general $1$-dimensional decomposition, one has solely
$\rho_\Gc(x, S_0)\geq n$,  for $x\in S_n$.
\end{example}

We define the \emph{inner boundary} of $S_n$ by
\[S_n^{-}:=\{x\in S_n,\, \exists y\in S_{n-1}, \,\Ec(x,y)>0\}\]
and the \emph{outer boundary} by
\[S_n^{+}:=\{x\in S_n,\, \exists y\in S_{n+1},\, \Ec(x,y)>0\}.\]

We now divide the degree with respect to the $1$-dimensional
decomposition. Given $x\in S_n$, we set
\begin{align*}
\deg_\pm(x):=\frac{1}{m(x)} \sum_{y\in S_{n \pm 1}} \Ec(x,y), \quad
\deg_0(x):=\frac{1}{m(x)} \sum_{y\in S_{n}} \Ec(x,y),
\end{align*}
with the convention that $S_{-1}=\emptyset$. Note that $\deg=
\deg_++\deg_-+\deg_0$ is independent of the choice of
a $1$-dimensional decomposition.

We now give a criterion of essential self-adjointness.

\begin{theorem}\label{t:1d}
Let $\Gc= (\Ec, \Vc, m)$ be a weighted graph and $(S_n)_{n\geq 0}$ be
a $1$-dimensional decomposition. Assume that
\[\sum_{n=1}^\infty \frac{1}{\sqrt{a_n^++ a_{n+1}^-} }=\infty,\]
where $a_n^\pm:= \sup_{x\in S_n^\pm}\deg_\pm(x)$.
Then $\Gc$ is $\chi$-complete and in particular, $\Delta_1$ is
essentially self-adjoint on $\Cc^c_{\rm   anti}(\Ec)$.
\end{theorem}
\proof
It is enough to check the hypothesis of Theorem \ref{t:25}. We set
\[\chi_n(x):=\left\{ \begin{array}{rl}
1, & \mbox{ if } |x|\leq n,
\\ \displaystyle
\max\left(0, 1- \sum_{k=n}^{|x|-1} \frac{1}{\sqrt{a_{k}^++ a_{k+1}^-}}\right),&
  \mbox{ if } |x|> n.
\end{array} \right.\]
Since the series diverges, $\chi_n$ is with
compact support. Note that $\chi_n$ is constant on $S_n$. If $x\in
S_m$ with $m>n$, we have~:
\begin{align*}
\frac{1}{m(x)}\sum_{y\in S_{m+1}} \Ec(x,y)|\chi_n(x)- \chi_n(y)|^2&\leq \deg_+(x)
\frac{1}{a_m^++ a_{m+1}^-}.
\\
\frac{1}{m(x)}\sum_{y\in S_{m}} \Ec(x,y)|\chi_n(x)- \chi_n(y)|^2&=0
\\
\frac{1}{m(x)}\sum_{y\in S_{m-1}} \Ec(x,y)|\chi_n(x)- \chi_n(y)|^2&\leq \deg_-(x)
\frac{1}{a_{m-1}^++ a_{m}^-}.
\end{align*}
It satisfies the definition of $\chi$-completeness. \qed

\subsubsection{Application to the class of antitrees}

We focus on antitrees, see
also~\cite{Wo2} and~\cite{BK}.
The \emph{sphere} of radius~$n\in\N$ around a
vertex~$v\in \Vc$ is the set $S_n(v):=\{w\in \Vc, \rho_\Gc(v,w)=n\}$.
A graph is an \emph{antitree}, if there exists a vertex $v\in \Vc$
such that for all other vertices $w\in \Vc\setminus\{v\}$
\begin{equation*}
  \Nc_\Gc(w)=S_{n-1}(v)\cup S_{n+1}(v)\text,
\end{equation*}
where $n=\rho_\Gc(v,w)\geq 1$.  See Figure \ref{fig_antitree} for an example.
The distinguished vertex~$v$ is the \emph{root} of the antitree.
Antitrees are bipartite and enjoy radial symmetry.

\begin{figure}
  \def\random{no}              % random sizes of spheres
  \def\f{1.5}                   % strech factor in x-direction
  \def\sizes{1, 7, 13, 5, 9, 20, 3}% (maximal) sizes of spheres
  \newcounter{n}                % n counts the spheres, starting at 0
  \begin{tikzpicture}
    \foreach \j in \sizes
    { \ifthenelse{\equal{\random}{yes}}{
        \pgfmathparse{random(\j)}
        \let\j\pgfmathresult
      }{}
      \ifthenelse{\j>1}{%
        \pgfmathparse{1-\j/2}\let\A\pgfmathresult
        \pgfmathparse{\j/2}\let\B\pgfmathresult
      }{\def\A{0.5}\def\B{1}}%
      \foreach \Y in {\A,...,\B} {
        \coordinate (x) at ({\f*(\then+1)},{(\Y-1/2)/sqrt(\B-\A)}) {};
        \filldraw (x) circle (0.7mm);
        \ifthenelse{\then=0}{}{%
          \foreach \y in {\a,...,\b}
            \draw (\f*\then,{(\y-1/2)/sqrt(\b-\a)}) -- (x);
        }%
      }
      \node [above=1.5mm] at (x) {$S_{\then}$}; % prints S_n
      \stepcounter{n} % next step
      \global\let\a\A % save coordinates for next loop
      \global\let\b\B % save coordinates for next loop
    }
  \end{tikzpicture}
  \addtocounter{n}{-1}
  \ifthenelse{\equal{\random}{yes}}
    {\caption{An antitree with spheres $S_0,\dotsc,S_{\then}$.}}
    {\caption{An antitree with spheres $S_0,\dotsc,S_{\then}$ of sizes $\sizes$.}}
  \label{fig_antitree}
\end{figure}

\begin{proposition}\label{p:anti}
Let $\Gc:(\Ec, \Vc)$ be a simple antitree whose root is $v$. Set
$s_n:=\sharp S_n(v)$.  Assume that
\[\sum_{n\in\N}
\frac{1}{\sqrt{s_{n}+s_{n+1}}}=\infty
\]
then $\Gc$ is $\chi$-complete and in particular, $\Delta_1$ is
essentially self-adjoint on $\Cc^c_{\rm   anti}(\Ec)$.
\end{proposition}
\proof Use Theorem \ref{t:1d} with $S_n:=S_n(v)$. \qed

\subsection{Other approaches}\label{s:other}
 We present some techniques that ensure essential self-adjointness.

\begin{theorem}\label{t:diff}
Let $\Gc=(\mathcal{E},\mathcal{V},m)$ be a locally finite graph. Set
\begin{align*}
\mathcal{M}(x,y):=1+\sum_{z\in\mathcal{V}}\left(\frac{1}{m(x)}\mathcal{E}(x,z)+\frac{1}{m(y)}\mathcal{E}(z,y)\right).
\end{align*}
 Suppose that
\begin{equation*}
\displaystyle
  \sup_{x,y\in\mathcal{V}}\sum_{z\in\mathcal{V}}\frac{1}{m(x)}\mathcal{E}(x,z)|
  \mathcal{M}(x,y)-\mathcal{M}(x,z) |^{2}<\infty.
\end{equation*}
Then $\Delta_{1}$ is essentially self-adjoint on $\mathcal{C}_{\rm anti}^{c}(\Ec)$.
\end{theorem}
\proof
Take $f\in \mathcal{C}_{\rm anti}^{c}(\Ec)$.  We have
\begin{align*} \| \Delta_{1}f \|^{2}\displaystyle
&\displaystyle\leq
\sum_{(x,y)\in\mathcal{V}^{2}}\mathcal{E}(x,y)\left( \frac{1}{m^{2}(x)}\Big|
\sum_{z\in\mathcal{V}}\mathcal{E}(x,z)f(x,z)\Big|^{2}\right.
\\
&\quad \quad \quad\quad\quad\quad\quad\quad\quad\quad\quad\quad\quad\quad\quad\left. +\frac{1}{m^{2}(y)}\Big|\sum_{z\in\mathcal{V}}\mathcal{E}(z,y)f(z,y)
\Big |^{2}\right)
\\
&=  2\sum_{(x,y)\in\mathcal{V}^{2}}\mathcal{E}(x,y)\frac{1}{m^{2}(x)}\Big |
\sum_{z\in\mathcal{V}}\mathcal{E}(x,z)f(x,z)\Big|^{2}
\\
&
\leq
2\sum_{(x,y)\in\mathcal{V}^{2}}\mathcal{E}(x,y)\frac{1}{m^{2}(x)}\Big(\sum_{t\in\mathcal{V}}\mathcal{E}(x,t)\Big)\Big(\sum_{z\in\mathcal{V}}\mathcal{E}(x,z)|
f(x,z)|^{2}\Big)
\\
&\displaystyle
=2\sum_{(x,y)\in\mathcal{V}^{2}}\mathcal{E}(x,y)\left|
\frac{1}{m(x)}\Big(\sum_{z\in\mathcal{V}}\mathcal{E}(x,z)\Big) f(x,y)\right|^{2}
\leq 2\|\mathcal{M}(\cdot, \cdot) f \|^2.
\end{align*}
  Moreover, noting that $\Mc(x,y)=\Mc(y,x)$ and that $f(x,y)=-f(y,x)$
 we get:
  \begin{align*}
  2|\langle f,[\Delta_{1},\mathcal{M}(\cdot, \cdot)]f\rangle|
&\displaystyle\leq
  \|f\|^2+
\\
&\hspace*{-1cm}\sum_{(x,y)\in\mathcal{V}^{2}}\mathcal{E}(x,y)\Big(\sum_{z\in\mathcal{V}}\frac{1}{m(x)}\mathcal{E}(x,z)\big|
\mathcal{M}(x,z)-\mathcal{M}(x,y) \big| \big|
  f(x,z) \big|\Big)^{2}
  \\
&\hspace{-1 cm} \leq \displaystyle\| \mathcal{M}(\cdot, \cdot)^{\frac{1}{2}}f\|^{2}
%\\
%&\hspace{-1 cm}
+\sum_{(x,y)\in\mathcal{V}^{2}}\mathcal{E}(x,y)\frac{1}{m(x)}\Big(\sum_{t\in\mathcal{V}}
\mathcal{E}(x,t)\Big) \times
\\
&\hspace{0cm}\Big(\sum_{z\in\mathcal{V}}\mathcal{E}(x,z)\frac{1}{m(x)}\big|\mathcal{M}(x,z)-\mathcal{M}(x,y)\big|^{2}\big|
f(x,z) \big|^{2}\Big)
\\
&\hspace{-1 cm} = \displaystyle\| \mathcal{M}(\cdot, \cdot)^{\frac{1}{2}}f\|^{2}
%\\
%&\hspace{-3.5 cm}
+
\sum_{(x,y)\in\mathcal{V}^{2}}\mathcal{E}(x,y)\sum_{t\in\mathcal{V}}\frac{1}{m(x)}\mathcal{E}(x,t)
\\
&
\underbrace{\sum_{z\in\mathcal{V}}\mathcal{E}(x,z)\frac{1}{m(x)}\big|\mathcal{M}(x,y)-\mathcal{M}(x,z)\big|^{2}}_{\leq C}
\big|f(x,y)\big|^{2}
\\
&\leq (1+2C)\| \mathcal{M}(\cdot, \cdot)^{\frac{1}{2}}f\|^{2}.\end{align*}
Applying \cite[Theorem X.37]{RS}, the result follows.
\qed

\begin{example} Take $\Gc=(\Ec, \Vc)$ to be a simple radial tree such
  that
\[\sup_{n\in \N}|\off(n)-\off(n+2)|<\infty, \]
 then $\Delta_1$ is
  essentially self-adjoint on $\Cc_{\rm anti}^c(\Ec)$.
 \end{example}

This example can be reached by Theorem \ref{t:25}, however for
weighted graphs the hypotheses do not seem to fully overlap from one
side or from the other.

We now give a more powerful way to check that $\Delta_{1}$ is
essentially self-adjoint. It seems that this approach is new in the
setting of graphs.
\begin{theorem}\label{t:n2}
Let $\Gc=(\Ec, \Vc, m)$ be a weighted graph. Fix $x_0\in \Vc$ and
set $B_n:=\{x\in \Vc, \rho_\Gc (x_0,x)\leq n\}$. If
\begin{align}\label{e:n2}
\sum_{n=1}^\infty \left( \prod_{i\leq n} \sup_{x\in B_i}(\deg_\Gc(x))
\right)^{-1/2n}= \infty
\end{align}
then $\Delta_1$ is essentially self-adjoint on $\Cc_{\rm
  anti}^c(\Ec)$.
\end{theorem}
\proof
Since $\Delta_1$ is non-negative, by \cite{MaMc} or \cite{Nu}, it is
enough to prove that for all $f\in \Cc_{\rm anti}^c(\Ec)$:
\[\sum_{n=0}^\infty \|\Delta_1^n f\|^{-1/2n}=\infty.\]
Set $f\in \Cc_{\rm anti}^c(\Ec)$. There is $n_0$ such that $\supp
f\subset B_{n_0}^2$, where $B_{n_0}^2:=B_{n_0}\times B_{n_0}$. Note
that $\supp (\Delta_1)^nf \subset
B_{n_0+n}^2$. Therefore
\begin{align*}
\|(\Delta_1)^{n+1} f\|& \leq \|1_{B_{n_0+n+1}^2}\Delta_1
1_{B_{n_0+n+1}^2}\|\cdot \|(\Delta_1)^{n} f\|
\\
&\leq \sup_{x\in
  B_{n_0+n+1}}\deg_\Gc(x)\cdot \|(\Delta_1)^{n} f\|
\end{align*}
since, for all $g\in \Cc_{\rm anti}^c(\Ec)$ and $k\in \N$, we have:
\begin{align*}
0\leq \langle 1_{B_{k}^2}g,\Delta_{1}1_{B_{k}^2}g\rangle
&=\displaystyle\sum_{x}\frac{1}{m(x)}\Big|\sum_{y}\mathcal{E}(x,y)
1_{B_{k}^2}(x,y)g(x,y)\Big|^{2}
\\
&\hspace{-2cm}\leq
\sum_{x}\frac{1}{m(x)}\Big(\sum_{y}\mathcal{E}(x,y)\Big)\sum_{z}
\mathcal{E}(x,z)|1_{B_{k}^2}(x,z)g(x,z)|^{2}
\\&\hspace{-2cm} \leq  \sup_{x\in B_k}\deg(x)\cdot\|1_{B_{k}^2} g \|^{2}.
\end{align*}
This concludes the proof. \qed

\begin{remark} The example of Proposition \ref{p:25sharp} is also covered by
  Theorem \ref{t:n2}.
\end{remark}

\begin{remark} A computation shows that the condition
  \eqref{e:n2} is equivalent to
\begin{align*}
\sum_{n=1}^\infty \frac{1}{\sqrt{\sup_{x\in B_n}(\deg_\Gc(x))}}= \infty,
\end{align*}
when $\sup_{x\in B_n}(\deg_\Gc(x))$ is equivalent to
$n^\alpha\ln^\beta(n)$, as $n$ goes to infinity.
\end{remark}

\subsection{The case of the Adjacency matrix}\label{s:adj}
We take the opportunity to apply this technique for the adjacency
matrix. This improves the results obtained in \cite{9}.
\begin{theorem}\label{t:n3}
Let $\Gc=(\Ec, \Vc, m)$ be a weighted graph. Fix $x_0\in \Vc$ and
set $B_n:=\{x\in \Vc, \rho_\Gc (x_0,x)\leq n\}$. If
\begin{align}\label{e:n3}
\sum_{n=1}^\infty \left( \prod_{i\leq n} \sup_{x\in B_i}(\deg_\Gc(x))
\right)^{-1/n}= \infty
\end{align}
then $\Ac_\Gc$ is essentially self-adjoint on $\Cc^c(\Vc)$.
\end{theorem}
\proof
Since $\Ac_{\Gc}$ is not necessary non-negative, by \cite{MaMc} or
\cite{Nu}, it is
enough to prove that for all $f\in \Cc^c(\Vc)$
\[\sum_{n=0}^\infty \|\Ac_{\Gc}^n f\|^{-1/n}=\infty.\]
Set $f\in \Cc^c(\Vc)$. There is $n_0$ such that $\supp
f\subset B_{n_0}$. Note that $\supp \Ac^n_{\Gc}f \subset
B_{n_0+n}$. Therefore
\begin{align*}
\|\Ac_{\Gc}^{n+1} f\|& \leq \|1_{B_{n_0+n+1}}\Ac_{\Gc}
1_{B_{n_0+n+1}}\|\cdot \|\Ac_{\Gc}^{n} f\|
\\
&\leq \sup_{x\in
  B_{n_0+n+1}}\deg_\Gc(x)\cdot \|\Ac^{n}_{\Gc} f\|
\end{align*}
since, for all $g\in \Cc^c(\Vc)$, we have
$|\langle g, \Ac_{\Gc} g\rangle|\leq \langle g, \deg_\Gc(\cdot) g\rangle$.\qed

\begin{corollary}\label{c:n3}
Let $\Gc=(\Ec, \Vc)$ be a radial simple tree,
endowed with an origin such that $\off(n):=\off(x)$ for
all $x\in \Sc_{n}$. Assume that
\[\sum_{n=1}^\infty \frac{1}{\off(n)}=\infty.\]
Then $\Ac_\Gc$ is essentially self-adjoint.
\end{corollary}

\begin{remark}
Note that the condition \eqref{e:n3} is not optimal for the question of
self-adjointness in the case of a radial simple tree, see
\cite{9}[Proposition 1.2] but is optimal in the
case of anti-trees, see \cite{GoSc2}.
\end{remark}

To finish denote by $\tilde \Gc$ the line graph of a
bipartite graph $\Gc$. Then, thanks to Proposition \ref{p:uni}, all the
criteria of Sections \ref{s:complete} and \ref{s:other} implies that
 $\Ac_{\tilde \Gc}$ is essentially self-adjoint on
$\Cc^c(\tilde \Vc)$.

\end{document}